\documentclass[a4paper]{article}
\usepackage[style=numeric]{biblatex}
\usepackage{graphicx}%
\usepackage{multirow}%
\usepackage{amsmath,amssymb,amsfonts}%
\usepackage{amsthm}%
\usepackage{authblk}%
\usepackage{mathrsfs}%
\usepackage[title]{appendix}%
\usepackage{xcolor}%
\usepackage{textcomp}%
\usepackage{manyfoot}%
\usepackage{booktabs}%
\usepackage{algorithm}%
\usepackage{algorithmicx}%
\usepackage{algpseudocode}%
\usepackage{listings}%
\usepackage{hyperref}

\usepackage{subcaption}%

\theoremstyle{thmstyleone}%

\begin{document}

\title{Grid-constrained online scheduling of flexible electric vehicle charging}

\author[1]{Emily van Huffelen}

\author[1]{Roel Brouwer}
% ORCID: 0000-0001-8300-4290

\author[1]{Marjan van den Akker}
% ORCID:

\affil[1]{Department of Information and Computing Sciences, Utrecht University, Utrecht, The Netherlands}

\begin{abstract}We study Electric Vehicle (EV) charging from a scheduling perspective, aiming to minimize delays while respecting the grid constraints. A network of parking lots is considered, each with a given number of charging stations for electric vehicles. Some of the parking lots have a roof with solar panels. The demand that can be served at each parking lot is limited by the capacity of the cables connecting them to the grid. We assume that EVs arrive at the parking lots according to a known distribution. Upon arrival, we learn the desired departure time, the amount of electrical energy it needs to charge its battery, and the range of rates that it can be charged at. Vehicle arrival patterns, connection times, and charging volume are based on data collected in the city of Utrecht. The departure time of an EV is delayed if it has not finished charging in time for its desired departure. We aim to minimize the total delay.
We present a novel approach, based on an online variant of well-known schedule generation schemes. We extend these schemes and include them in a destroy-and-repair heuristic. This resulted in several scheduling strategies. We show their effectiveness using a discrete event simulation. With this, we show that applying scheduling approaches increases the amount of EVs that can be charged at a site and reduces the average delay. Furthermore, we argue the importance of considering aspects of the grid layout in electricity networks and show the benefits of using flexible charging rates.
\end{abstract}

\maketitle

\section{Introduction}
The number of electric vehicles (EVs) has rapidly increased in recent years. This trend is likely to continue, as policies aiming to reduce greenhouse gas emissions encourage their sale to replace vehicles that run on fossil fuels. The increasing number of EVs requires charging stations to be installed in many locations. At parking lots where many EVs congregate, such as one near a large office building for example, managing the demand for charging the batteries of all these EVs concurrently is a challenge. Such parking lots often have clear peaks in demand, and are relatively quiet at other times. For parking lots outside of residential neighborhoods, this peak in demand often occurs during the day, which matches well with the peak in production of solar panels. Therefore, the installation of a large number of charging stations is often paired with the installation of solar panel arrays on roofs over these parking lots. Although the peaks in supply and demand do not line up as nicely for parking lots in residential areas, solar panels are commonly installed in these areas as well. The servicing of EVs at these parking lots will benefit from smart scheduling of the charging jobs to align well with the availability of power at the parking lot.

In this work, we consider a network of parking lots operated by the same provider, who offers the service of parking and charging electric vehicles (EVs).
Each parking lot has a given number of charging stations for EVs. The demand that can be served at each parking lot is limited by the amount of available power. The power supply is constrained by the capacity of the cables connecting the parking lots to the grid. Some parking lots are covered by a roof with solar panels on it, supplying additional power under the right weather conditions.

We assume that EVs arrive at the parking lots according to a known distribution. As soon as a vehicle arrives, we learn its desired departure time, the amount of electrical energy it needs to charge its battery before that time, and the range of rates that it can be charged at. The charging of these EVs must be scheduled in such a way that the total delay is minimized, while the network constraints are respected. This means that the summed difference between the desired departure time and the actual departure time, which is delayed if an EV has not finished charging in time, must be minimal across all vehicles.

The charging of electric vehicles is a widely studied subject across multiple disciplines, even if we limit ourselves to the (optimal) scheduling of charging jobs. Many surveys aim to provide an overview of aspects of this problem. Wang et al. \cite{Wang2016} make a classification of EV charging control algorithms based on their perspective and objective. They identify three categories: \textit{smart grid oriented}, \textit{aggregator oriented} and \textit{customer oriented} EV charging, describing the perspective from which the problem is considered. Although the objectives discussed in the review are mainly cost optimization objectives, we could consider our problem setting to align with the \textit{aggregator oriented} model, as we aim to provide the best possible service within the limited capacity of the network.

The control of EV charging can be considered on several `levels', or time scales \cite{Arif2021}. The lowest level of control is reactive, and mainly concerns stabilization of the grid. Arif et al. \cite{Arif2021} describe that EVs can be used at that level to provide ancillary services by controlling the charging (grid-to-vehicle, G2V) or discharging (vehicle-to-grid, V2G) at a very small time scale. Such time scales make scheduling intractable. We focus on higher levels of control: aligning demand well with forecasts of solar panel output and the demand of other EVs within the same network. The importance of the uncertainty of demand is emphasized by Al-Ogaili et al. \cite{AlOgaili2019}. This implies that offline scheduling methods may struggle to account for the uncertainty surrounding the arrival times and charging demands of EVs at charging locations. Therefore, we consider an online scheduling approach.

Most surveys conclude that coordinated charging strategies, where some centralized decision making is involved, yield the best results. Hussain et al. \cite{Hussain2021}, for example, states that ``it is found that centralized coordination is best strategy to handle all issues effectively." Liu et al. \cite{Liu2015} also find that centralized approaches have better performance, but prefer decentralized approaches due to the communication overhead of centralized approaches.

We will be looking at a centralized control strategy, from the perspective of an aggregator. We aim to schedule the demand in a network of parking places as well as possible, given the limited resources available. In this sense, the problem we consider is closely related to the Resource-Constrainted Project Scheduling Problem (RCPSP). The surveys by Hartmann and Briskorn \cite{Hartmann2010,Hartmann2022} give a good overview of the RCPSP and its variants. The RCPSP with flexible resource profiles (FRCPSP) such as studied by Naber \cite{Naber2017} and the General Continuous Energy-Constrained Scheduling Problem (GCECSP) \cite{Brouwer2024}

The approaches we develop in this work are inspired by priority rule based scheduling heuristics commonly employed for finding solutions to RCPSP and its variants. Kolisch \cite{Kolisch1996} provides a good overview of the application of schedule generation schemes to the classical RCPSP. Lova et al. \cite{Lova2006} apply schedule generation schemes to the multi-mode RCPSP, where tasks can be processed in a number of different modes, which can be seen as a discretized variant of the flexible model we are considering. More recently, priority rule based approaches that are similar to schedule generation schemes have been applied to the stochastic RCPSP \cite{Chen2018} and in the context of genetic algorithms, where compound priority rules are learned for the generation of schedules \cite{Dumic2021}.

\textbf{Our contribution:} We present a novel approach, based on an extension of traditional schedule generation schemes to an online setting with flexible charging rates. We consider single-pass methods, where we apply a scheme once, as well as more advanced methods, where we include the schemes in a destroy-and-repair heuristic. In this way, we develop a number of variants of the approach for the generation of efficient schedules for charging EVs on a network of parking lots. As far as we are aware, similar approaches have not been applied to this problem yet. We show the effectiveness of these approaches using a discrete event simulation. Furthermore, we argue the importance of considering aspects of the grid layout in electricity networks and show the benefits of using flexible charging rates. 

The work presented here is a continuation of earlier work by van Huffelen \cite{vanHuffelen2023}. The specific case is loosely based on an assignment used in the master's course `Optimization for sustainability' at Utrecht University \cite{vandenAkker2021}. The data used for the distribution of arrival times, charging volumes and connection times is based on real-world data from the city of Utrecht, and was provided by Brinkel et al. \cite{Brinkel2020}. The distribution of solar output is based on a forecast of the production of solar panels in the Netherlands in 2025 \cite{ENTSOE2018}.

The rest of this work is structured as follows. In Section \ref{sec:evcharging-prob-desc}, we give a detailed problem description, followed by an explanation of the developed scheduling approaches in Section \ref{sec:evcharging-scheduling}. We present the results of evaluating our approaches using that simulation in Section \ref{sec:evcharging-results}. Finally, we will draw some conclusions and discuss avenues for future work in Section \ref{sec:evcharging-conclusion}.

\section{Problem description}\label{sec:evcharging-prob-desc}

\begin{figure}
	\centering
	\includegraphics[]{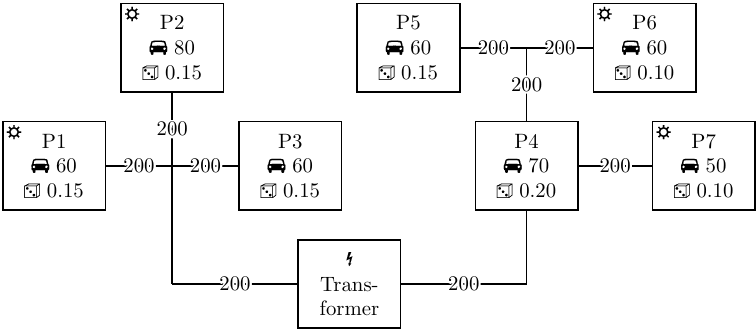}
	\caption{Network layout}
	\label{fig:evcharging-network}
\end{figure}

We consider an online problem, where a schedule for the charging of EVs needs to be generated with the objective to minimize the average delay. EVs arrive over time, which implies that the schedule needs to be updated regularly. We will discuss strategies for these updates in Section \ref{sec:evcharging-scheduling}.

We consider a network that consists of a number of interconnected parking lots, that share a grid connection. An example network is shown in Figure \ref{fig:evcharging-network}. Each parking lot has:
\begin{itemize}
	\item a number of parking spots, where EVs can be parked and charged;
	\item the fraction of EVs $p_{\text{P}i}$ that will select parking lot P$i$;
	\item (optionally) an array of solar panels, producing a variable amount of power.
\end{itemize}

These parking lots are connected by cables with a given maximum capacity. The cable capacity limits the amount of power that can be transported to each parking lot. We use this a simple approximation of the limits imposed by the network.

The amount of power generated by the solar panels is based on the so-called availability vector. This expresses the average production of solar panels at any given time as a fraction of their peak power.
The production of the solar panels is revealed at the start of every hour, and treated as a constant production rate for the entire hour.

The EVs arrive at this location according to a Poisson process. The arrival rate is not constant, but changes every hour. Each EV $j$ has a number of properties, which are revealed when the vehicle arrives. Here, EV $j$ is the $j$th vehicle in the order of arrival:
\begin{itemize}
	\item a connection time $d_j - r_j$, the amount of time between its arrival time $r_j$ and its desired time of departure $d_j$;
	\item a charging volume $E_j$: the amount of energy required to fully charge its battery before departure (in kWh);
	\item charging rates $[P^-_j,P^+_j]$, a range of rates between a lower ($P^-_j$) and an upper ($P^+_j$) bound that can be used to charge it, in kW;
	\item a parking preference: a ordered list of three parking lots, where the EV will try to park. If all three parking lots in the parking preference of an EV are full, it will leave the area.
\end{itemize}
For each EV, we determine a start and completion time for the charging process. The start time must be after its arrival, which implies that the EV may not start charging immediately. After its start, we determine a charging profile that will be followed. It completes when its full energy demand is served, which will be earlier than or, if that is not possible, as close to its preferred departure time as possible.
Note that we guarantee that the connection time is never smaller than the time it takes to charge an EV at its minimum rate, i.e. $E_j/P^-_j$.

In addition, note that the existence of a lower bound implies that preemption is not allowed: as soon as an EV has started charging, it cannot drop below the lower bound on its charging rate before it completes. A reason for this is that interruption of the charging process comes with some loss of efficiency caused by startup and shutdown effects.

In this setting, we want to minimize the total tardiness $\sum_j \max(0, C_j - d_j)$ over all EVs that parked at the network of parking lots.

In this work, we use a case study for the evaluation of our approaches, loosely based on an assignment used in the master's course `Optimization for sustainability' at Utrecht University \cite{vandenAkker2021}. It uses the network depicted in Figure \ref{fig:evcharging-network}. All cables have a capacity of 200 kW. We use the availability factor predicted for solar panels in summer in the Netherlands in 2025 \cite{ENTSOE2018}, which we provide in Appendix \ref{app:evcharging-distributions-solar}. We assume the average hourly power output to follow a normal distribution, with the value from the availability vector as its mean, and a standard deviation of 15\%. All solar panel arrays installed have a peak power of 200 kW, unless explicitly stated otherwise.
We assume an average arrival rate of 1125 EVs per day. We use hourly arrival rates, charging volumes, and connection times based on data from Brinkel et al. \cite{Brinkel2020}. These distributions are provided in Appendix \ref{app:evcharging-distributions-arrival},  \ref{app:evcharging-distributions-volumes} and \ref{app:evcharging-distributions-connection}. The distribution of charging ranges can be found in Appendix \ref{app:evcharging-distributions-rates}.

\section{Scheduling}\label{sec:evcharging-scheduling}
As EVs keep arriving continuously, we use an online scheduling approach. Only the characteristics of the EVs that have already arrived have been revealed, and the amount of power produced by the solar panels is only known for the current hour. The available information will serve as the basis for a schedule that determines when and at what rate each EV charges. In this section, we will first introduce the methods that are used to find schedules based on this information. After that, we will discuss a strategy for replacing these schedules when new information is revealed.

\subsection{Schedule generation schemes}\label{sec:evcharging-sgs}
In this chapter, we consider priority rule based scheduling heuristics. These are made up of two components: a schedule generation scheme and a priority rule. Priority rules will be discussed in Section \ref{sec:evcharging-priority-rules}. The schedule generation scheme is a constructive heuristic that adds jobs to the schedule one by one, until a feasible schedule is obtained.

There are two main variants of schedule generation schemes that are often considered in literature: serial and parallel \cite{Kolisch1996}. A convenient way to characterize the difference is that the serial scheme iterates over jobs, while the parallel scheme iterates over points in time. The serial scheme aims to schedule the highest priority job as early as possible, and the parallel scheme tries to fill the current time with as high-priority jobs as possible. Both variants construct a schedule by adding jobs, until all jobs have been placed. In simple schedule generation schemes, the schedule that results from this procedure is not adapted further. These schemes are considered to be \textit{single-pass} schedule generation schemes, since they place all jobs in a single loop over jobs (serial) or time (parallel).

At any stage in the construction, the serial scheme will select the unscheduled job with highest priority (according to the priority rule) and schedule it as early as possible. It will then repeat this until all jobs are scheduled. The parallel scheme, in contrast, will select the highest priority job (or jobs, if enough power remains available after the scheduling of the first one) that fits in the gap at the current time in the schedule, which is not necessarily the highest priority job available overall.

We extended these schemes in the following ways to apply them to our case:
\begin{description}
	\item[Composite resource] To check whether enough resource (power) is available to process a job at a given rate, we need to check the load of all cables and parking lots along the path from its source to the lot where the EV is located. The residual capacity at a parking lot can be determined by an elementary flow computation from the load at other parking lots, including the contribution of solar panels, and the cable capacities.
	\item[Resilience against drops of solar power] We need to guarantee the load on the cables does not exceed the capacity, while using as much of the available power as we can. If, for whatever reason, the solar output should fall away, we want to be able to guarantee that we can scale back charging on some EVs to remain within the bounds of the cable capacity. Therefore, we ensure that the EV's that are currently charging, can charge at their minimum rate without using solar power, i.e. solar power is only used to increase the charging rate above the lower bound.
	\item[Flexible rate] The jobs have a flexible resource consumption profile. This means that we need to decide the consumption rate for each job for the duration of its execution. When a job is selected to be added, it is placed at the earliest possible time in the partial schedule, at the highest possible rate. It may start at any rate between its lower and upper bound. From that point, it will always be scheduled to charge at the highest possible rate, until it completes. This rate is equal to its upper bound, or the maximum rate allowed by the resource constraints, whichever is more restrictive. This means the rate may change over time.
	\item[Avoid preemption] Applying the schemes in an online setting means that we need to take into account jobs that are already started when we reschedule. Recall that preemption is not allowed. When the construction of a new schedule starts, these jobs are initialized as charging at their lowest possible rate. As soon as they are encountered in the priority order, this rate is scaled up where possible. In this way, preemption is avoided, while respecting the priority order as much as possible.
\end{description}

Note that in the parallel schedule generation scheme, we do not need to generate a complete schedule, but just one that lasts until the next time a schedule will be generated in our strategy. We stop the construction as soon as it advances to a time point beyond that. This is not true for the serial scheme, as a lower priority job may still cause changes early in the schedule.

%The performance of these adaptations of the traditional single-pass serial and parallel schedule generation schemes will be evaluated in Section \ref{sec:evcharging-results}.

\subsection{Priority rules}\label{sec:evcharging-priority-rules}
A schedule generation scheme uses a priority rule to assign a priority to each job. A priority rule essentially is a formula that assigns a numerical value to a job, that can be used to order them. Typically, jobs with a smaller value for the priority rule's formula have a higher priority. Many possible priority rules exist, but we will only list the ones that will be used in the remainder of this work:
\begin{description}
	\item[FCFS] First Come First Serve: $r_j$.
	\item[EDD] Earliest Due Date: $d_j$.
	\item[ELSTu] Updated Earliest Latest Starting Time: $d_j - E'_j/P^+_j$.
\end{description}
where $E'_j$ is the remaining charging volume, i.e. the original charging volume $E_j$ minus the amount charged up to the time of evaluation.

\subsection{Destroy-and-repair}\label{sec:evcharging-scheduling-dar}
The single-pass schedule generation schemes described in Section \ref{sec:evcharging-sgs} generate reasonably good schedules. We want to investigate, however, if we can improve them further by using an iterative heuristic.

Essentially, our approach boils down to a \textit{destroy-and-repair} heuristic. We remove parts of the schedule and replace them by parts that are generated using a different priority rule. We repeat this a limited number of times to obtain an improved schedule.

\begin{algorithm}[H]
	\caption{Global overview of the destroy-and-repair procedure. $J$ is the set of available jobs (of size $n$), $t$ is the current time.}\label{alg:evcharging-rescheduling}
	\begin{algorithmic}[1]
		\State $\mathcal{S}^* \gets $ InitialSchedule($J$, PRIO1) \Comment{Generate initial schedule with a parallel or serial schedule generation scheme using the first priority rule}
		\State fails $\gets 0$
		\While{$\text{fails} < f$}
			\State $\mathcal{S} \gets \mathcal{S}^*$
			\State $J' \gets \emptyset$
			\For{$\text{iter} \in \{1, ..., r \cdot n\}$} \Comment{Randomly remove $r \cdot n$ jobs}
				\State $j \gets $ GetRandomJob($\mathcal{S}$)
				\State $J' \gets J'\cup\{j\}$
				\State $\mathcal{S} \gets $ RemoveFromSchedule($\mathcal{S}, j$)
			\EndFor
			\For{$\text{iter} \in \{r \cdot n, ..., s \cdot n\}$} \Comment{Further removal based on adjacency}
				\State $j \gets $ GetWeightedRandomJob($\mathcal{S}$)
				\State $J' \gets J'\cup\{j\}$
				\State $\mathcal{S} \gets $ RemoveFromSchedule($\mathcal{S}, j$)
			\EndFor
			\For{$j \in J' : S_j < t$} \Comment{Reinsert already started jobs at minimum rate}
				\State $\mathcal{S} \gets $ ScheduleAtMinRate($\mathcal{S}, j$)
			\EndFor
			\For{$j \in$ SortBy($J'$, PRIO2)} \Comment{Reinsert (or adjust) jobs in order of the second priority rule}
				\State $\mathcal{S} \gets $ AddToSchedule($\mathcal{S}, j$)
			\EndFor
			\If{Score($\mathcal{S}$) $<$ Score($\mathcal{S}^*$) + $i_{\text{min}}$} \Comment{Above some minimum improvement threshold $i_{min}$, we count a success}
				\State fails $ \gets$ fails + 1
			\Else
				\State fails $ \gets 0$
			\EndIf
			\If{Score($\mathcal{S}$) $>$ Score($\mathcal{S}^*$)} \Comment{Keep the best result}
				\State $\mathcal{S}^* \gets \mathcal{S}$
			\EndIf
		\EndWhile
		\Return{$\mathcal{S}^*$}
	\end{algorithmic}
\end{algorithm}

A global overview of the procedure is given in Algorithm \ref{alg:evcharging-rescheduling}. From the initial schedule, we remove a fraction $s$ of jobs. First a fraction $r < s$ is removed randomly with uniform probability. For the remaining fraction $s-r$, jobs are removed with a probability proportional to their weight. The weight is defined as the number of jobs adjacent to the job in the original schedule that have already been removed. In this way, the removal procedure aims to create some large contiguous open areas in the schedule, rather than many small holes. For two jobs to be considered adjacent, two conditions must be met: (1) part of their processing windows must overlap or immediately follow each other, and  (2) they must be scheduled on two parking lots that compete for resources (i.e. parking lots that share at least one cable along the path to the root node).

We start the repair phase by reinserting jobs that have a start time $S_j$ before the start of the schedule, i.e. jobs that are already being processed, to avoid preemption. These are initially scheduled at the lowest possible rate. Then, we reinsert (or adjust) all jobs, ordering them using a different priority rule than the one used to generate the initial schedule.

The resulting schedule is compared to the best schedule we have seen so far. If it is better, we continue the procedure with the new schedule. If it is better by at least a minimum improvement $i_{\text{min}}$, we say the procedure was successful. If the procedure was not successful a number of $f$ consecutive times, we stop.

\subsection{Scheduling frequency}
We consider an online problem. This means that new information regularly becomes available, either when the solar output is updated, or when a new EV arrives at a parking lot. The most naive way to deal with this is to generate a new schedule any time new information is revealed. However, since we use a rather advanced scheduling algorithm, it is reasonable to generate a new schedule less often. Therefore, we trigger the generation of a new schedule on a periodic basis, e.g. every hour. Within an interval, no new schedules are generated, unless a high-priority EV arrives. A high-priority EV is defined as an EV that has a higher priority value (for the first priority rule) than 80\% of the EVs in the schedule at that time. If such an EV arrives, an additional schedule is generated at its time of arrival. It does not change the timing of the next periodic schedule generation, however.

\section{Computational results}\label{sec:evcharging-results}
%Gesimuleerde periode in de code: 9 dagen - 2 dagen warm-up
We have implemented a discrete-event simulation to evaluate our approach. A brief description of the simulation is presented in Appendix \ref{app:evcharging-events}.

For the test runs reported below, average results over five runs are presented, each with a fixed seed used to generate the instance. Each individual run contains nine days of simulated time, the first two of which are treated as a warm-up period. This means that all reported results (except the runtime) are over days 3-9. The simulator was written in the Python programming language. The processor of the system used to run the tests on was an Intel(R) Xeon(R) Gold 6130 CPU @ 2.10GHz.

The main performance metric of an approach is the resulting average delay, or the average tardiness, of an EV that is charged at one of the parking lots in the simulation. The other metrics that we will report on are maximum delay, i.e. the largest delay a single EV encounters in the simulation, the fraction of EVs that is delayed as a percentage of the total number of EVs that parked at one of the parking lots, and the number of EVs that suffered a delay of more than fifteen minutes. We present the number of EVs with a large delay rather than a percentage, as we feel that the \emph{number} of EVs that experience a large delay is a relevant metric, even if the percentage is small.

In the remainder of this section, we first discuss the selected approaches and parameter settings in \ref{sec:evcharging-parameters}. We continue with an analysis of the effects of considering elements of the grid topology (Section \ref{sec:evcharging-results-parking}) and flexible charging rates (Section \ref{sec:evcharging-results-rate}). Then we will show the performance of the developed approaches on different scheduling intervals in Section \ref{sec:evcharging-results-interval}. We conclude this section with a comparison of the most promising approaches in Section \ref{sec:evcharging-results-general}.

\subsection{Approach selection and parameter settings}\label{sec:evcharging-parameters}
% Iets zeggen over zon, arrivals, ...
Based on their performance in initial test runs, we selected five variants that we will use in our evaluation:
\begin{description}
	\item[FCFS] is a serial schedule generation scheme, using the \textit{First Come First Serve} priority rule.
	\item[S] generates schedules using a serial schedule generation scheme, where the planned departure of an e-vehicle is used as the priority (Earliest Due Date/EDD).
	\item[P] generates schedules using a parallel schedule generation scheme, where the planned departure of an e-vehicle is used as the priority (Earliest Due Date/EDD).
	\item[SR] generates initial schedules using a serial schedule generation scheme, where the planned departure of an e-vehicle is used as the priority (Earliest Due Date/EDD). Improvements are generated using the destroy-and-repair heuristic described in Section \ref{sec:evcharging-scheduling-dar} with the Updated Earliest Latest Starting Time (ELSTu) priority rule.
	\item[PR] generates initial schedules using a parallel schedule generation scheme, where the planned departure of an e-vehicle is used as the priority (Earliest Due Date/EDD). Improvements are generated using the destroy-and-repair heuristic described in Section \ref{sec:evcharging-scheduling-dar} with the Updated Earliest Latest Starting Time (ELSTu) priority rule.
\end{description}
For the latter two, the values of the relevant parameters are as follows:
\begin{itemize}
	\item Rescheduling fraction $s = 0.5$;
	\item Random removal fraction $r = 0.05$;
	\item Improvement threshold $i_{\text{min}} = 100$;
	\item Unsuccessful iteration tolerance $f = 4$.
\end{itemize}

Other variants of the approaches have been tested, using different (combinations) of priority rules, parameter settings, and adjacency measures (in case of weighted removal). The presented approaches and settings (except `FCFS') were chosen as they showed the most consistent good performance in the preliminary tests. Apart from these four main approaches, we also present results from `FCFS' when relevant.

\subsection{Network constraints}\label{sec:evcharging-results-parking}
\begin{figure}
	\centering
	\includegraphics{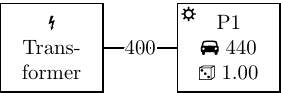}
	\caption{Copperplate layout}
	\label{fig:evcharging-copperplate}
\end{figure}
% Methods: P-EDD S-EDD PR-EDD-LSTu SR-EDD-LSTu

Where limitations exist in the transport capacity of the underlying electricity network, it is important to take these into account when planning electricity usage. To show the effects of ignoring the grid layout, we have implemented a scenario that does so (see Figure \ref{fig:evcharging-copperplate}). This network consists of a single parking lot, with as many parking spots and solar panels as the seven parking lots in our original grid combined. The results of running our four main approaches on both grids are shown in Table \ref{tab:evcharging-results-grid}.

We observe that delays almost disappear when we take away the grid layout. This indicates that the simplification of ignoring the grid layout gives a significant and unrealistic advantage in the scheduling model. This difference justifies the additional modeling effort required to account for the topology of the network.

\begin{table}
	\centering
		\begin{tabular}{ll|rrrr}
			\toprule
			& Metric & S & P & SR & PR \\
			\midrule
			\multirow{4}{*}{Grid} & Max. delay (s) & 14236.27 & 13882.98 & 5033.37 & 4847.46 \\
			& Avg. delay (s) & 121.19 & 106.03 & 88.14 & 86.56 \\
			& EVs delayed (\%) & 1.75 & 1.48 & 1.47 & 1.54 \\
			& Delay $\geq$ 15 min & 65.20 & 58.00 & 52.60 & 55.40 \\
			\midrule
			\multirow{4}{*}{No Grid} & Max. delay (s) & 36.41 & 0.00 & 36.41 & 10.40 \\
			& Avg. delay (s) & 0.01 & 0.00 & 0.01 & 0.00 \\
			& EVs delayed (\%) & 0.00 & 0.00 & 0.00 & 0.00 \\
			& Delay $\geq$ 15 min & 0.00 & 0.00 & 0.00 & 0.00 \\
			\bottomrule
		\end{tabular}
		\caption{Comparison of results taking the network layout into account and those ignoring the grid}
		\label{tab:evcharging-results-grid}
\end{table}

\subsection{Flexible charging rates}\label{sec:evcharging-results-rate}
The assumption that using the flexibility in rates that charging of EV batteries may allow would be beneficial, is crucial for our model. To test this assumption, we compared simulations with the flexible rates as we have described them so far, to simulations with a fixed rate of 9 kW for all EVs. The value of 9 kW was chosen, as this is the average upper bound in the flexible scenario, meaning that the average potential for energy consumption is equal across both scenarios.
Additionally, we excluded the solar panels from this simulation. EVs cannot scale down in the fixed rate scenario and their charging is not allowed to be preempted. Because of this, using the power from solar panels in the schedule means that the cable capacity may be violated if the amount of power they provide changes. To make the comparison as fair as possible, therefore, the solar panels are turned off. Note the absence of solar power causes markedly different results compared to e.g. Table \ref{tab:evcharging-results-grid}, as less power is available for charging overall. The results of these runs are presented in Table \ref{tab:evcharging-results-flexible}.

Here, we present the results of the `FCFS' alongside our four main approaches, as it performs unexpectedly well in this scenario. It seems that planning according to the `first come first serve' principle is fairly effective in a strained network where all power flows from a single source.

We observe that the differences are not large, but the average delay is typically a few minutes shorter in the scenario with flexible charging rates. The same is true for the amount of severe delays and the percentage of EVs that is delayed.

Keep in mind that we excluded the solar panels from this comparison. Their output can be incorporated safely in the scenario with flexible rates, whereas it cannot be incorporated when using fixed charging rates without allowing preemption. Therefore, one could say that flexible rates are very beneficial when we want to use solar power. These results together show that the use of flexible charging rates allows for better outcomes.

\begin{table}
	\centering
		\begin{tabular}{@{}l@{ \ }l|r@{ \ }r@{ \ }r@{ \ }r@{ \ }r@{}}
			\toprule
			& Metric & FCFS & S & P & SR & PR \\
			\midrule
			\multirow{4}{*}{Fixed} & Max. delay (s) & 68148.35 & 54599.73 & 54561.94 & 53366.20 & 56653.89 \\
			& Avg. delay (s) & 15142.88 & 14899.12 & 14894.50 & 15055.51 & 15103.98 \\
			& EVs delayed (\%) & 61.30 & 85.83 & 85.89 & 85.83 & 85.01 \\
			& Delay $\geq$ 15 min & 2106.40 & 2950.80 & 2954.40 & 2958.80 & 2911.20 \\
			\midrule
			\multirow{4}{*}{Flexible} & Max. delay (s) & 66068.74 & 62563.74 & 59144.77 & 61743.67 & 60528.36 \\
			& Avg. delay (s) & 14131.81 & 14600.97 & 14486.47 & 14327.40 & 14418.19 \\
			& EVs delayed (\%) & 59.02 & 82.34 & 83.25 & 81.61 & 83.86 \\
			& Delay $\geq$ 15 min & 2055.60 & 2837.40 & 2870.60 & 2815.60 & 2886.60 \\
			\bottomrule
		\end{tabular}
		\caption{Comparison of results for flexible and fixed rate charging}
		\label{tab:evcharging-results-flexible}
\end{table}

\subsection{Scheduling frequency}\label{sec:evcharging-results-interval}
In an online setting, regular updates of the schedule are performed. Each time new information is revealed, the generation of a new schedule could be triggered. The results presented in the previous sections use this strategy, because it provides the most clear environment for the comparison of elements of the modeling. In practice, however, it is undesirable and unnecessary to generate new schedules that often, in particular for the more advanced scheduling methods. Generating schedules less frequently is also more realistic in terms of computation time, as new EVs arrive with high frequency during peak hours.

The more often a new schedule is generated, the better the results are expected to be. If the difference in quality is sufficiently small, however, it is reasonable to update schedules with a lower frequency. We have tested three scheduling frequencies:
\begin{enumerate}
	\item Generate a new schedule whenever new information is revealed (i.e., when the solar forecast is updated or a new EV arrives);
	\item Generate a new schedule once every 15 minutes;
	\item Generate a new schedule once every hour.
\end{enumerate}
In the latter two cases, we keep track of the priority of newly arriving EVs, according to the (first) priority rule used for scheduling. If an EV arrives with a priority that is higher than 80\% of the EVs that are currently scheduled, we will trigger an additional schedule to be generated.

The results are presented in Table \ref{tab:evcharging-results-interval}. Most approaches show a slight deterioration in the quality of the results when the scheduling frequency decreases. Only the results for the parallel scheduling scheme (`P') are somewhat inconsistent with this observation. The performance of the `SR' approach deteriorates most of all.

We observe that the total runtime of the simulation is significantly reduced for all approaches. Note that the reported runtime includes every aspect of the simulation (not just the scheduling itself) for a total simulation time of nine days. The runtime of all approaches is such that it is realistic to implement them with a schedule generation interval of 15 minutes. In terms of the quality of the schedules, we observe that the best performing (advanced) approach (`PR') still outperforms the single-pass approaches on the main metric (average delay), even with a large scheduling interval.
\begin{table}
	\centering
		\begin{tabular}{ll|rrrr}
			\toprule
			& Metric & S & P & SR & PR \\
			\midrule
			\multirow{5}{*}{New information} & Max. delay (s) & 14236.27 & 13882.98 & 5033.37 & 4847.46 \\
			& Avg. delay (s) & 121.19 & 106.03 & 88.14 & 86.56 \\
			& EVs delayed (\%) & 1.75 & 1.48 & 1.47 & 1.54 \\
			& Delay $\geq$ 15 min & 65.20 & 58.00 & 52.60 & 55.40 \\
			& Runtime (s) & 20366.03 & 12743.83 & 70506.54 & 40210.12 \\
			\midrule
			\multirow{5}{*}{15 min.} & Max. delay (s) & 14706.39 & 13167.38 & 5293.81 & 5237.77 \\
			& Avg. delay (s) & 123.94 & 115.11 & 106.39 & 95.23 \\
			& EVs delayed (\%) & 1.75 & 1.63 & 1.73 & 1.60 \\
			& Delay $\geq$ 15 min & 65.20 & 58.00 & 52.60 & 55.40 \\
			& Runtime (s) & 8688.63 & 7143.22 & 27071.67 & 15324.67 \\
			\midrule
			\multirow{5}{*}{1 hour} & Max. delay (s) & 14889.79 & 14027.18 & 8399.59 & 5135.09 \\
			& Avg. delay (s) & 127.69 & 103.22 & 136.22 & 96.95 \\
			& EVs delayed (\%) & 1.77 & 1.47 & 1.94 & 1.61 \\
			& Delay $\geq$ 15 min & 65.20 & 58.00 & 52.60 & 55.40 \\
			& Runtime (s) & 7499.83 & 7423.24 & 19654.80 & 14446.14 \\
			\bottomrule
		\end{tabular}
		\caption{Comparison of results for a number of scheduling intervals.}
		\label{tab:evcharging-results-interval}
\end{table}

\subsection{Overall performance}\label{sec:evcharging-results-general}
From the results in the previous sections, we observe that the parallel implementation of an approach seems to dominate the serial implementation of that same approach. Based on the results of our preliminary tests, however, we must say that this is not a general result. This is consistent with earlier findings by e.g. Kolisch \cite{Kolisch1996}.

In this section we will give a final overview of the results. We present the best performing single-pass approach and the best performing approach using our destroy-and-repair heuristic. We compare them to the `FCFS' approach, which we consider to be the most straight-forward application of scheduling using our adapted schedule generation schemes. The best performing approaches are the parallel schedule generation scheme using the Earliest Due Date (EDD) priority rule, `P-EDD', and the parallel `improved' approach using the EDD and ELSTu (Updated Earliest Latest Starting Time) priority rules, `PR-15min'. For the latter approach, we use a scheduling interval of 15 minutes, while the other two will continuously generate new schedules, any time new information is revealed.

The parallel `improved' approach outperforms the other approaches on almost all metrics, while it is comparable in computational load. This seems to be the most promising approach for our case study and similar scenarios. 

Note that all approaches discussed above protect the network capacity and apply scheduling techniques for charging decisions. We conclude with a note on the effectiveness of scheduling in general. Therefore, we want to make a comparison with a scenario where no control is exercised over the charging at all.
We ran simulations on the same network, under the same conditions, where all EVs have a fixed charging rate of 9 kW. Rather than scheduling their demand, any EV that arrives at a parking lot immediately starts charging, until its battery is fully charged. In this scenario, we consider a 10\% overload of the cable capacity (i.e. a load of 200-220 kW) to be manageable. We consider anything beyond 220 kW to be truly problematic. If we look at the two main cables connecting the two groups of parking lots to the transformer, we find that the cable on the left (connecting to P1, P2 and P3) suffers from a manageable overload 4.4\% of the time, and from a problematic overload 53.7\% of the time. The cable of the right (connecting to P4, P5, P6 and P7) has a manageable overload 2.6\% of the time, and a problematic overload 65.8\% of the time. These numbers indicate that it is absolutely necessary to manage the charging of EVs on these parking lots, and show the usefulness and effectiveness of a scheduling approach like ours.

\begin{table}
	\centering
		\begin{tabular}{l|rrr}
			\toprule
			& FCFS & P-EDD & PR-15min \\
			\midrule
			Max. delay (s) & 31287.96 & 13882.98 & 5237.77 \\
			Avg. delay (s) & 1114.21 & 106.03 & 95.23 \\
			EVs delayed (\%) & 11.72 & 1.48 & 1.60 \\
			%EVs served (\%) & 54.54 & 55.48 & 55.55 \\
			Delay $\geq$ 15 min & 470.20 & 58.00 & 55.40 \\
			Runtime (s) & 16168.70 & 12743.83 & 15324.67 \\
			\bottomrule
		\end{tabular}
		\caption{Results for the most promising approaches}
		\label{tab:evcharging-results-full}
\end{table}

\section{Conclusions and future work}\label{sec:evcharging-conclusion}
We have formulated an online scheduling problem dealing with the charging of EVs. For this problem, we developed a scheduling approach that extends the idea of schedule generation schemes to an online setting with flexible charging jobs. In our case study, the application of schedules generated using this approach, allows for the serving of a large number of EVs with minimal delay, while respecting the constraints imposed by the network. If no control over the charging was exercised at all, cables in the network would be overloaded around 60-70\% of the time. Overload is completely avoided if our schedules are followed, at the cost of an average delay of just over 1.5 minutes per EV in a high-occupancy scenario. The quality of the generated schedules is such that frequent updates are not necessary. Finally, the computational load of the developed approach is reasonable for practical purposes, making it a promising candidate for implementation in real-world cases.

In addition to this main result, our case study shows that it is essential to take the network topology into account when deciding charging profiles. Furthermore, we have shown the benefits of considering a range of permissible charging rates for EVs, rather than a single fixed value.

Many avenues for future research exist. Broadly, we see two main directions for future research.

First, their are many opportunities for further development of the scheduling approach. We have proposed a novel approach for generating schedules. Other variants of this approach might be interesting. Examples include the use of more than two priority rules in the improvement process and the application of a more directed destroy-and-repair heuristic that aims to improve the schedule at the points that cause the most delay.

Second, certain extensions of the model may be interesting to investigate. For example, a more complex model could include the ability of EVs to discharge as well (V2G), and contribute to the charging of other EVs with a higher priority. A more detailed modeling of the uncertain power output of the solar panels can be considered as well. Using a more detailed model of the power network (e.g. an AC model) would be an interesting, but computationally heavy, extension. Other objectives, such as cost minimization through consideration of electricity prices, can also be investigated.

\section*{Acknowledgements}

We would like to thank Philip de Bruin for providing the numbers for the simulation runs without any control or scheduling strategies, that we used in Section \ref{sec:evcharging-results-general}.
We want to acknowledge the research program ``Energie: Systeem Integratie en Big Data? with project number 647.003.005, which is financed by the Dutch Research Council (NWO), for financing the machine used for the computational experiments.

\section*{Declarations}

\subsection*{Funding}
No funding was received to assist with the preparation of this manuscript.

\subsection*{Competing interests}
The authors have no competing interests to declare that are relevant to the content of this article.

\subsection*{Data availability}
The data used in this paper is provided in the form of tables in the appendix.

\subsection*{Code availability}
The code used to generate the results presented in this paper is not available publicly, as parts of it may be used for assignments in a course. Access to the code can be requested from the corresponding author.

\subsection*{Author contributions}
All authors contributed to the study conception and design. Coding and analysis were performed by Emily van Huffelen under supervision of Marjan van den Akker and Roel Brouwer. Further analysis was done by Roel Brouwer, who also wrote the first draft of the manuscript. All authors commented on previous versions of the manuscript. All authors read and approved the final manuscript.

\begin{appendices}

\section{Simulation model}\label{app:evcharging-events}
We evaluated the developed scheduling approaches using a discrete-event simulation. We will briefly explain the simulation model using the event graph displayed in Figure \ref{fig:evcharging-event-graph}. Each box represents a type of event, where arcs between events indicate that an event may schedule another. A dashed arc indicates that the subsequent event happens instantaneously if it is triggered. The arcs labeled with `(a)' may or may not be present, depending on the scheduling approach used (see Section \ref{sec:evcharging-scheduling}). The simulation is run by repeatedly extracting the next event from the event queue and executing its associated actions, until the \textit{end simulation} event is encountered.

We simulate the arrival, charging and departure of electric vehicles. The charging is done according to the schedule we generate during the simulation, using one of the approaches described in Section \ref{sec:evcharging-scheduling}. The shaded boxes in Figure \ref{fig:evcharging-event-graph} indicate events that are contained in the schedule. A schedule is, in essence, a series of provisional events of types \textit{EV changes charging rate} and \textit{EV stops charging}. The \textit{inspect schedule} event draws these events from the schedule and enters them in the event queue. A brief description of each event is provided below:

\begin{figure}[h]
	\centering
	\includegraphics{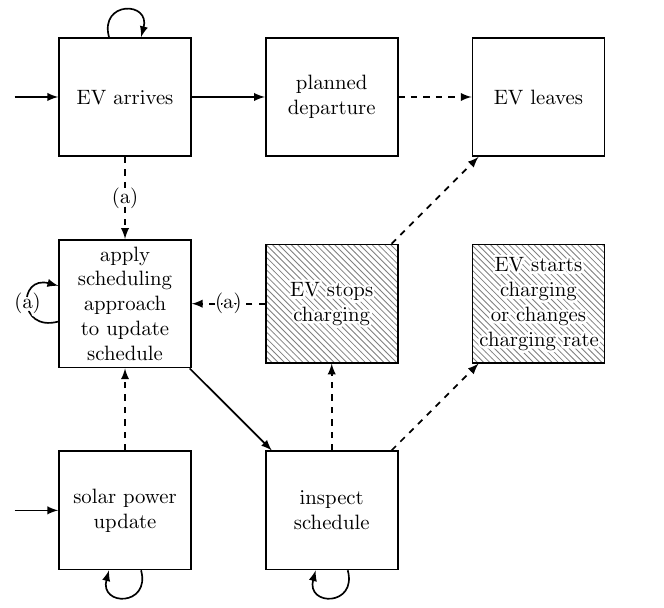}
	\caption{Event graph showing the relations between events in our simulation model}
	\label{fig:evcharging-event-graph}
\end{figure}

\begin{description}
	\item[EV arrives] happens when a new EV enters the simulation. This event handles the selection of a parking lot (or removes the EV from the simulation if all three parking lots from its preference are full). It schedules a \textit{planned departure} event for the EV at its due date $d_j$ and may schedule an immediate \textit{update schedule} event, depending on the scheduling approach used. Finally, it will schedule a \textit{EV arrives} event for the next EV to arrive.
	\item[planned departure] happens when an EV is supposed to leave, at the end of its planned connection time. If the EV is fully charged, it will schedule an immediate \textit{EV leaves} event.
	\item[EV leaves] happens when an EV actually leaves. This means that it leaves its parking lot and is removed from the simulation.
	\item[EV starts charging or changes charging rate] happens when an EV starts charging for the first time, or changes its charging rate later on. It updates all relevant values and rates, but schedules no additional events. This event is dictated by the schedule.
	\item[EV stops charging] happens when an EV stops charging, i.e. it has received its requested charging volume $E_j$. It updates all relevant values and rates. If the due date $d_j$ of the EV has already passed, it will schedule an immediate \textit{car leaves} event. Depending on the scheduling approach used, it may also trigger an immediate \textit{update schedule} event.
	\item[solar power update] is a regularly scheduled event, that happens on a given interval: every hour of simulated time. This event retrieves the (expected) power produced by the solar panels installed for the upcoming interval, and updates the available power at each parking lot accordingly. It will trigger an immediate \textit{update schedule} event. Finally, it will schedule another \textit{regular forecast update} event at the end of the current forecast interval.
	\item[apply scheduling approach to update schedule] happens when the schedule needs to be updated. A new schedule is generated using the relevant scheduling approach. It schedules an \textit{inspect schedule} event at the time of occurrence of the first action in the generated schedule. Additionally, if schedules are updated on an interval different from the \textit{regular forecast update}, it schedules another \textit{update schedule} event at the end of the current scheduling interval.
	\item[inspect schedule] will retrieve the next event from the schedule (either a \textit{EV changes charging rate} event or a \textit{EV stops charging} event) and, if it happens at the current time, schedules the corresponding event. It will then schedule an \textit{inspect schedule} event at the (planned) time of occurrence of the following event in the schedule.
\end{description}

At the start of the simulation, three events are scheduled:
\begin{enumerate}
	\item `regular forecast update' at time 0;
	\item `car arrives' at the first generated arrival time;
	\item `end simulation' at the end of the period of time we intend to simulate.
\end{enumerate}
These correspond to the three arcs that do not originate at an event box in Figure \ref{fig:evcharging-event-graph}.
\section{Distributions}\label{app:evcharging-distributions}
\subsection{Solar power}\label{app:evcharging-distributions-solar}
The power that a solar panel produces strongly depends on the time of day. In Table \ref{tab:evcharging-distributions-solar}, we list the average production of solar panels as a fraction of their peak power. These will be used to determine the mean of the distribution (with a standard deviation of 15\%) from which the actual values are drawn during the simulation. The values are based on a forecast of the production of solar panels during the summer months in the Netherlands in 2025 \cite{ENTSOE2018}.

\begin{table}[hp]
	\centering
	\begin{minipage}[b]{0.3\linewidth}
	\begin{tabular}{rr}
		\toprule
		Hour & Fraction \\
		\midrule
		0 & 0.000000000 \\
		1 & 0.000000000 \\
		2 & 0.000000000 \\
		3 & 0.001548952 \\
		4 & 0.017126799 \\
		5 & 0.055567520 \\
		6 & 0.132034210 \\
		7 & 0.222931250 \\
		\bottomrule
	\end{tabular}
	\end{minipage}
	\quad
	\begin{minipage}[b]{0.3\linewidth}
	\begin{tabular}{rr}
		\toprule
		Hour & Fraction \\
		\midrule
		8 & 0.311153080 \\
		9 & 0.382618250 \\
		10 & 0.426692930 \\
		11 & 0.446001000 \\
		12 & 0.440104200 \\
		13 & 0.403637380 \\
		14 & 0.342447800 \\
		15 & 0.261458840 \\
		\bottomrule
	\end{tabular}
	\end{minipage}
	\quad
	\begin{minipage}[b]{0.3\linewidth}
	\begin{tabular}{rr}
		\toprule
		Hour & Fraction \\
		\midrule
		16 & 0.169520680 \\
		17 & 0.083306720 \\
		18 & 0.026936987 \\
		19 & 0.005260382 \\
		20 & 0.000000000 \\
		21 & 0.000000000 \\
		22 & 0.000000000 \\
		23 & 0.000000000 \\
		\bottomrule
	\end{tabular}
	\end{minipage}
	\caption{Average solar panel revenues as a fraction of peak power} \label{tab:evcharging-distributions-solar}
\end{table}

\subsection{Arrival of EVs}\label{app:evcharging-distributions-arrival}
EVs arrive in the simulation according to a Poisson process. The (average) arrival rates for every hour of the day as a fraction of the total daily arrivals are listed in Table \ref{tab:evcharging-distributions-arrival}.
\begin{table}[hp]
	\centering
	\begin{minipage}[b]{0.3\linewidth}
	\begin{tabular}{rr}
		\toprule
		Hour & Fraction \\
		\midrule
		0 & 0.012926435 \\
		1 & 0.004938874 \\
		2 & 0.002042621 \\
		3 & 0.001188988 \\
		4 & 0.000304869 \\
		5 & 0.000518277 \\
		6 & 0.007560745 \\
		7 & 0.027072345 \\
		\bottomrule
	\end{tabular}
	\end{minipage}
	\quad
	\begin{minipage}[b]{0.3\linewidth}
	\begin{tabular}{rr}
		\toprule
		Hour & Fraction \\
		\midrule
		8 & 0.058748209 \\
		9 & 0.049114356 \\
		10 & 0.035456236 \\
		11 & 0.040120728 \\
		12 & 0.045181549 \\
		13 & 0.048108289 \\
		14 & 0.050913082 \\
		15 & 0.057254352 \\
		\bottomrule
	\end{tabular}
	\end{minipage}
	\quad
	\begin{minipage}[b]{0.3\linewidth}
	\begin{tabular}{rr}
		\toprule
		Hour & Fraction \\
		\midrule
		16 & 0.072223408 \\
		17 & 0.103899271 \\
		18 & 0.125453492 \\
		19 & 0.071156367 \\
		20 & 0.056217798 \\
		21 & 0.052467913 \\
		22 & 0.047315631 \\
		23 & 0.029816164 \\
		\bottomrule
	\end{tabular}
	\end{minipage}
	\caption{Average fraction of arrivals at each hour} \label{tab:evcharging-distributions-arrival}
\end{table}
\subsection{Charging rates}\label{app:evcharging-distributions-rates}
Each EV has a minimum and maximum charging rate (kW), drawn with uniform probability from the combinations listed in Table \ref{tab:evcharging-distributions-rates}.

\begin{table}[hp]
	\centering
	\begin{minipage}[b]{0.45\linewidth}
	\begin{tabular}{rrr}
		\toprule
		Min. & Max. & Probability \\
		\midrule
		3 & 6 & 0.0625 \\
		3 & 7 & 0.0625 \\
		3 & 8 & 0.0625 \\
		3 & 9 & 0.0625 \\
		4 & 7 & 0.0625 \\
		4 & 8 & 0.0625 \\
		4 & 9 & 0.0625 \\
		4 & 10 & 0.0625 \\
		\bottomrule
	\end{tabular}
	\end{minipage}
	\quad
	\begin{minipage}[b]{0.45\linewidth}
	\begin{tabular}{rrr}
		\toprule
		Min. & Max. & Probability \\
		\midrule
		5 & 8 & 0.0625 \\
		5 & 9 & 0.0625 \\
		5 & 10 & 0.0625 \\
		5 & 11 & 0.0625 \\
		6 & 9 & 0.0625 \\
		6 & 10 & 0.0625 \\
		6 & 11 & 0.0625 \\
		6 & 12 & 0.0625 \\
		\bottomrule
	\end{tabular}
	\end{minipage}
	\caption{Distribution of charging rates}
	\label{tab:evcharging-distributions-rates}
\end{table}
\subsection{Charging volumes}\label{app:evcharging-distributions-volumes}
Each EV has a demand for a given amount of energy to be added to its battery over the duration of its connection. These are distributed in the range of 0 to 102 kWh according to the fractions given in Table \ref{tab:evcharging-distributions-volumes}. These fractions represent the probability of the charging volume to be drawn from that particular range. Within each group, values are distributed uniformly.

\begin{table}[hp]
	\centering
	\begin{minipage}[b]{0.3\linewidth}
	\begin{tabular}{rr}
		\toprule
		kWh & Weight \\
		\midrule
		0-1 & 0.029938112 \\
		1-2 & 0.026188226 \\
		2-3 & 0.035212341 \\
		3-4 & 0.032681930 \\
		4-5 & 0.043443797 \\
		5-6 & 0.050333831 \\
		6-7 & 0.071186854 \\
		7-8 & 0.055973903 \\
		8-9 & 0.048809488 \\
		9-10 & 0.032163654 \\
		10-11 & 0.021188378 \\
		11-12 & 0.023383433 \\
		12-13 & 0.020670102 \\
		13-14 & 0.019725008 \\
		14-15 & 0.018688455 \\
		15-16 & 0.018261638 \\
		16-17 & 0.018048230 \\
		17-18 & 0.017347032 \\
		18-19 & 0.017286058 \\
		19-20 & 0.018718941 \\
		20-21 & 0.018078717 \\
		21-22 & 0.018901863 \\
		22-23 & 0.017286058 \\
		23-24 & 0.017316545 \\
		24-25 & 0.017103137 \\
		25-26 & 0.016615347 \\
		26-27 & 0.014511753 \\
		27-28 & 0.015487333 \\
		28-29 & 0.014938569 \\
		29-30 & 0.013200817 \\
		30-31 & 0.011859395 \\
		31-32 & 0.013322765 \\
		32-33 & 0.010853328 \\
		33-34 & 0.011127710 \\
		\bottomrule
	\end{tabular}
	\end{minipage}
	\quad
	\begin{minipage}[b]{0.3\linewidth}
	\begin{tabular}{rr}
		\toprule
		kWh & Weight \\
		\midrule
		34-35 & 0.009603366 \\
		35-36 & 0.008444864 \\
		36-37 & 0.008292430 \\
		37-38 & 0.008079022 \\
		38-39 & 0.006585165 \\
		39-40 & 0.007042468 \\
		40-41 & 0.006920521 \\
		41-42 & 0.006676626 \\
		42-43 & 0.006219323 \\
		43-44 & 0.005762019 \\
		44-45 & 0.005822993 \\
		45-46 & 0.005487638 \\
		46-47 & 0.005640072 \\
		47-48 & 0.005883967 \\
		48-49 & 0.005121795 \\
		49-50 & 0.004115728 \\
		50-51 & 0.004115728 \\
		51-52 & 0.004268163 \\
		52-53 & 0.003871833 \\
		53-54 & 0.004786439 \\
		54-55 & 0.003810859 \\
		55-56 & 0.003262096 \\
		56-57 & 0.003384043 \\
		57-58 & 0.003445017 \\
		58-59 & 0.003201122 \\
		59-60 & 0.002804793 \\
		60-61 & 0.002987714 \\
		61-62 & 0.002957227 \\
		62-63 & 0.002317003 \\
		63-64 & 0.002073108 \\
		64-65 & 0.002256029 \\
		65-66 & 0.001981647 \\
		66-67 & 0.001585318 \\
		67-68 & 0.001615804 \\
		\bottomrule
	\end{tabular}
	\end{minipage}
	\quad
	\begin{minipage}[b]{0.3\linewidth}
	\begin{tabular}{rr}
		\toprule
		kWh & Weight \\
		\midrule
		68-69 & 0.001249962 \\
		69-70 & 0.001128014 \\
		70-71 & 0.000640224 \\
		71-72 & 0.000853633 \\
		72-73 & 0.000518277 \\
		73-74 & 0.000731685 \\
		74-75 & 0.000518277 \\
		75-76 & 0.000762172 \\
		76-77 & 0.000579251 \\
		77-78 & 0.000396329 \\
		78-79 & 0.000396329 \\
		79-80 & 0.000182921 \\
		80-81 & 0.000182921 \\
		81-82 & 0.000091500 \\
		82-83 & 0.000091500 \\
		83-84 & 0.000091500 \\
		84-85 & 0.000274382 \\
		85-86 & 0.000243895 \\
		86-87 & 0.000152434 \\
		87-88 & 0.000182921 \\
		88-89 & 0.000182921 \\
		89-90 & 0.000182921 \\
		90-91 & 0.000030500 \\
		91-92 & 0.000000000 \\
		92-93 & 0.000000000 \\
		93-94 & 0.000030500 \\
		94-95 & 0.000000000 \\
		95-96 & 0.000000000 \\
		96-97 & 0.000000000 \\
		97-98 & 0.000000000 \\
		98-99 & 0.000000000 \\
		99-100 & 0.000000000 \\
		100-101 & 0.000000000 \\
		101-202 & 0.000030500 \\
		\bottomrule
	\end{tabular}
	\end{minipage}

	\caption{Distribution of charging volumes of EVs}\label{tab:evcharging-distributions-volumes}
\end{table}
\subsection{Connection times}\label{app:evcharging-distributions-connection}
Each EV has an intended time of departure, effectively a due date for the charging to be complete. This is determined by adding a connection time (in hours, distributed as in Table \ref{tab:evcharging-distributions-connection}) to the arrival time. These fractions represent the probability of the connection time to be drawn from that particular range. Within each group, values are distributed uniformly.
\begin{table}[hp]
	\centering
	\begin{minipage}[b]{0.3\linewidth}
	\begin{tabular}{rr}
		\toprule
		Hours & Weight \\
		\midrule
		0-1 & 0.075089174 \\
		1-2 & 0.089052163 \\
		2-3 & 0.084570592 \\
		3-4 & 0.069967379 \\
		4-5 & 0.053565440 \\
		5-6 & 0.036309869 \\
		6-7 & 0.025608975 \\
		7-8 & 0.028200360 \\
		8-9 & 0.041492637 \\
		9-10 & 0.037559830 \\
		10-11 & 0.037072040 \\
		11-12 & 0.040699979 \\
		12-13 & 0.046431511 \\
		13-14 & 0.055242218 \\
		14-15 & 0.049205817 \\
		15-16 & 0.034511143 \\
		16-17 & 0.026828450 \\
		17-18 & 0.022590775 \\
		18-19 & 0.019237218 \\
		19-20 & 0.016737295 \\
		20-21 & 0.013261791 \\
		21-22 & 0.010426511 \\
		22-23 & 0.008536325 \\
		23-24 & 0.006707113 \\
		\bottomrule
	\end{tabular}
	\end{minipage}
	\quad
	\begin{minipage}[b]{0.3\linewidth}
	\begin{tabular}{rr}
		\toprule
		Hours & Weight \\
		\midrule
		24-25 & 0.005518124 \\
		25-26 & 0.004237676 \\
		26-27 & 0.002926740 \\
		27-28 & 0.001859699 \\
		28-29 & 0.001524344 \\
		29-30 & 0.001128014 \\
		30-31 & 0.000579251 \\
		31-32 & 0.000914606 \\
		32-33 & 0.001280449 \\
		33-34 & 0.001249962 \\
		34-35 & 0.001341423 \\
		35-36 & 0.001920673 \\
		36-37 & 0.002256029 \\
		37-38 & 0.003231609 \\
		38-39 & 0.003536478 \\
		39-40 & 0.002804793 \\
		40-41 & 0.002408463 \\
		41-42 & 0.002865766 \\
		42-43 & 0.001890186 \\
		43-44 & 0.002286516 \\
		44-45 & 0.001890186 \\
		45-46 & 0.001310936 \\
		46-47 & 0.001432883 \\
		47-48 & 0.001341423 \\
		\bottomrule
	\end{tabular}
	\end{minipage}
	\quad
	\begin{minipage}[b]{0.3\linewidth}
	\begin{tabular}{rr}
		\toprule
		Hours & Weight \\
		\midrule
		48-49 & 0.000945093 \\
		49-50 & 0.000548764 \\
		50-51 & 0.000518277 \\
		51-52 & 0.000274382 \\
		52-53 & 0.000121948 \\
		53-54 & 0.000152434 \\
		54-55 & 0.000426816 \\
		55-56 & 0.000243895 \\
		56-57 & 0.000213408 \\
		57-58 & 0.000457303 \\
		58-59 & 0.000396329 \\
		59-60 & 0.000274382 \\
		60-61 & 0.000670711 \\
		61-62 & 0.000945093 \\
		62-63 & 0.000792659 \\
		63-64 & 0.000884119 \\
		64-65 & 0.000853633 \\
		65-66 & 0.000640224 \\
		66-67 & 0.000396329 \\
		67-68 & 0.000609738 \\
		68-69 & 0.000548764 \\
		69-70 & 0.000457303 \\
		70-71 & 0.007987561 \\
		&             \\
		\bottomrule
	\end{tabular}
	\end{minipage}

	\caption{Distribution of connection times} \label{tab:evcharging-distributions-connection}
\end{table}

\end{appendices}

\printbibliography

\end{document}